\documentclass[a4paper,12pt]{amsart}
\usepackage{amssymb}
\usepackage{ifthen}
\usepackage[usenames]{color}
\usepackage{graphicx}
\nonstopmode \numberwithin{equation}{section}
\newtheorem{thm}{Theorem}[section]

\newtheorem{cor}{Corollary}[section]
\newtheorem{lem}{Lemma}[section]
\newtheorem{prop}{Proposition}[section]

\newtheorem{claim}{Claim}[section]
\newtheorem{subclaim}{Subclaim}
\newtheorem{conj}[equation]{Conjecture}
\newtheorem{case}{Case}[section]
\newtheorem*{mysolution}{Solution}
\newtheorem{step}{Step}[section]
\theoremstyle{definition}
\newtheorem{defn}{Definition}[section]
\newtheorem{examp}{Example}[section]
\newtheorem{prob}[equation]{Problem}
\newtheorem{ques}[equation]{Question}
\newtheorem{rem}{Remark}[section]
\newcounter {own}
\def\theown {\thesection       .\arabic{own}}

\newenvironment{pf}[1][]{%
 \vskip 3mm
 \noindent
 \ifthenelse{\equal{#1}{}}%
  {{\slshape Proof. }}%
  {{\slshape #1.} }%
 }%
{\qed\bigskip}

\newcounter{alphabet}
\newcounter{tmp}
\newenvironment{Thm}[1][]{\refstepcounter{alphabet}%
\bigskip%
\noindent%
{\bf Theorem \Alph{alphabet}}%
\ifthenelse{\equal{#1}{}}{}{ (#1)}%
{\bf .} \itshape}{\vskip 8pt}

\makeatletter

\newcommand{\Ref}[1]{\@ifundefined{r@#1}{}{\setcounter{tmp}{\ref{#1}}\Alph{tmp}}}
\newenvironment{Lem}[1][]{\refstepcounter{alphabet}%
\bigskip%
\noindent%
{\bf Lemma \Alph{alphabet}}%
{\bf .} \itshape}{\vskip 8pt}

\newcommand{\IR}{{\mathbb R}}

\newcommand{\IB}{{\mathbb B}}

\def\be{\begin{equation}}
\def\ee{\end{equation}}

\newcommand{\ben}{\begin{enumerate}}
\newcommand{\een}{\end{enumerate}}

\newcommand{\blem}{\begin{lem}}
\newcommand{\elem}{\end{lem}}
\newcommand{\bthm}{\begin{thm}}
\newcommand{\ethm}{\end{thm}}
\newcommand{\bcor}{\begin{cor}}
\newcommand{\ecor}{\end{cor}}
\newcommand{\beg}{\begin{examp}}
\newcommand{\eeg}{\end{examp}}
\newcommand{\begs}{\begin{examples}}
\newcommand{\eegs}{\end{examples}}
\newcommand{\bdefe}{\begin{defn}}
\newcommand{\edefe}{\end{defn}}
\newcommand{\bprob}{\begin{prob}}
\newcommand{\eprob}{\end{prob}}
\newcommand{\bques}{\begin{ques}}
\newcommand{\eques}{\end{ques}}
\newcommand{\bei}{\begin{itemize}}
\newcommand{\eei}{\end{itemize}}
\newcommand{\bcl}{\begin{claim}}
\newcommand{\ecl}{\end{claim}}
\newcommand{\bscl}{\begin{subclaim}}
\newcommand{\escl}{\end{subclaim}}
\newcommand{\bca}{\begin{case}}
\newcommand{\eca}{\end{case}}
\newcommand{\bstep}{\begin{step}}
\newcommand{\estep}{\end{step}}
\newcommand{\bsol}{\begin{mysolution}}
\newcommand{\esol}{\end{mysolution}}
\newcommand{\bcon}{\begin{conj}}
\newcommand{\econ}{\end{conj}}
\newcommand{\bcons}{\begin{conjs}}
\newcommand{\econs}{\end{conjs}}
\newcommand{\bprop}{\begin{prop}}
\newcommand{\eprop}{\end{prop}}
\newcommand{\br}{\begin{rem}}
\newcommand{\er}{\end{rem}}
\newcommand{\brs}{\begin{rems}}
\newcommand{\ers}{\end{rems}}
\newcommand{\bo}{\begin{obser}}
\newcommand{\eo}{\end{obser}}
\newcommand{\bos}{\begin{obsers}}
\newcommand{\eos}{\end{obsers}}
\newcommand{\bpf}{\begin{pf}}
\newcommand{\epf}{\end{pf}}
\newcommand{\ba}{\begin{array}}
\newcommand{\ea}{\end{array}}
\newcommand{\beq}{\begin{eqnarray}}
\newcommand{\beqq}{\begin{eqnarray*}}
\newcommand{\eeq}{\end{eqnarray}}
\newcommand{\eeqq}{\end{eqnarray*}}

\begin{document}
\bibliographystyle{amsplain}

\title{Generalized Bloch spaces, Integral means of hyperbolic harmonic mappings in the unit ball}

\author{Jiaolong Chen}
\address{Jiaolong Chen, Key Laboratory of High Performance Computing and Stochastic Information Processing (HPCSIP)(Ministry of Education of China),
College of Mathematics and Computer Science,
Hunan Normal University, Changsha, Hunan 410081, People's Republic of China}
\email{jiaolongchen@sina.com}

%\author{Xiantao Wang${}^{~\mathbf{*}}$}
%\address{Xiantao Wang, Department of Mathematics,
%Shantou University, Shantou, Guangdong 515063, People's Republic of China}
%\email{xtwang@stu.edu.cn}

\subjclass[2000]{Primary: 31B05; Secondary: 31C05}
\keywords{Hyperbolic harmonic mapping, generalized Bloch space, integral means, weak uniform boundedness property.%\\
%$^{\mathbf{*}}$Corresponding author
}

%\dedicatory{}
\begin{abstract} In this paper, we investigate the properties of hyperbolic harmonic mappings in the unit ball $\mathbb{B}^{n}$ in $\IR^n$ $(n\geq 2)$.
Firstly, we establish necessary and sufficient conditions for a hyperbolic harmonic mapping to be in the Bloch space $\mathcal{B}(\mathbb{B}^{n})$ and the generalized Bloch space $\mathcal{L}_{\infty,\omega}\mathcal{B}_{\alpha,\mathrm{a}}^{0}(\mathbb{B}^{n})$, respectively.
Secondly, we discuss the relationship between the integral means of hyperbolic harmonic mappings and that of their gradients. The obtained results are the generalizations of Hardy and Littlewood's related ones in the setting of hyperbolic harmonic mappings.
Finally, we characterize the weak uniform boundedness property of hyperbolic harmonic mappings in terms of the quasihyperbolic metric.
\end{abstract}

\thanks{The research was partly supported by NSF of China (No. 11571216 and No. 11671127) and the construct program of the key discipline in Hunan Province}

\maketitle \pagestyle{myheadings} \markboth{Jiaolong Chen}{Generalized Bloch spaces, Integral means of hyperbolic harmonic mappings in the unit ball}

%%%%%%%%%%%%%%%%%%%%%%%%%%%%%%%%%%%%
%%%%%%%%%%%%%%%%%%%%%%%%%%%%%%%%%%%%
\section{Introduction and main results}\label{sec-1}
%%%%%%%%%%%%%%%%%%%%%%%%%%%%%%%%%%%%
%%%%%%%%%%%%%%%%%%%%%%%%%%%%%%%%%%%%

For $n\ge 2$, let $\mathbb{B}^{n}(x_{0}, r)=\{x\in\mathbb{R}^{n}:|x-x_{0}|<r\}$, $\mathbb{S}^{n-1}(x_{0}, r)=\partial\mathbb{B}^{n}(x_{0}, r)$ and $\overline{\mathbb{B}}^{n}(x_{0}, r)=\mathbb{B}^{n}(x_{0}, r)\cup \mathbb{S}^{n-1}(x_{0}, r)$. In particular, we write $\mathbb{B}^{n}=\mathbb{B}^{n}(0, 1)$, $\mathbb{S}^{n-1}=\mathbb{S}^{n-1}(0, 1)$ and $\overline{\mathbb{B}}^{n}=\mathbb{B}^{n}\cup \mathbb{S}^{n-1}$.

The purpose of this paper is to consider the hyperbolic harmonic mappings whose definition is as follows.

\bdefe\label{def-1.1}
A mapping $u=(u_1,\cdots,u_n)\in C^{2}(\mathbb{B}^{n}, \mathbb{R}^{n})$ is said to be {\it hyperbolic harmonic} if
 $$\Delta_{h}u=(\Delta_{h}u_{1}, \cdots,\Delta_{h}u_{n})=0,$$  that is,
 for each $j\in \{1,\cdots, n\}$, $u_j$ satisfies the hyperbolic Laplace equation
 $$\Delta_{h}u_{j}=0,$$ where
\be\label{eq-01.01}
\Delta_{h}u_{j} (x)=(1-|x|^2)^2\Delta u_{j}(x)+2(n-2)(1-|x|^2)\sum_{i=1}^{n} x_{i} \frac{\partial u_{j}}{\partial x_{i}}(x).
\ee
 \edefe
\noindent We refer to \cite{chrw, grja, jami, reka, sto1999, sto2012, sto2016} for basic properties of this class of mappings. For convenience, in the following of this paper, we always use the notation $
\Delta_{h}u=0$ to mean that $u=(u_1,\cdots,u_n)$ is hyperbolic harmonic in $\mathbb{B}^{n}$.

Obviously, for $n=2$, hyperbolic harmonic mappings coincide with harmonic mappings. See \cite{clsh,du} and the references therein for the basic properties of harmonic mappings.

%%%%%%%%%%%%%%%%%%%%%%%%%%%%%%%%%%%%
\subsection{Generalized Hardy spaces and generalized Bloch spaces}
%%%%%%%%%%%%%%%%%%%%%%%%%%%%%%%%%%%%

 \bdefe\label{def-1.2}
For $p\in (0,\infty]$, the {\it generalized Hardy space $\mathcal{H}_{g}^{p}(\mathbb{B}^{n} )$} consists of all those functions
$f: \mathbb{B}^{n}\rightarrow\mathbb{ R}^{n}$ such that each $f_{i}$ is measurable, $M_{p}(r,f)$ exists for all $r\in (0,1)$ and $$||f||_{p}<\infty,$$ where
$$f=(f_1,\cdots,f_n),\;\;\;\;\;||f||_{p}=\sup_{0<r<1} \big\{M_{p}(r,f)\big\}$$
and
$$\;\;M_{p}(r,f)=\begin{cases}
\displaystyle \;\left( \int_{\mathbb{S}^{n-1}} |f(r\xi)|^{p}d\sigma(\xi)\right)^{\frac{1}{p}} , & \text{ if } p\in (0,\infty),\\
\displaystyle \;\sup_{\xi\in \mathbb{S}^{n-1}} \big\{|f(r\xi)|\big\} , \;\;\;\;& \text{ if } p=\infty.
\end{cases}$$
\edefe
\noindent Here and hereafter, $d \sigma$ always denotes the normalized surface measure on $\mathbb{S}^{n-1}$ so that $\sigma(\mathbb{S}^{n-1})=1$.

The classical {\it Hardy space $\mathcal{H}^{p}(\mathbb{D} )$} consisting of related analytic functions is a subspace of
$\mathcal{H}_{g}^{p}(\mathbb{D} )$, where $\mathbb{D}$ denotes the unit disk in the complex plane $\mathbb{C}$ (In this paper, we always identify $\mathbb{R}^2$ with $\mathbb{C}$ and $\mathbb{B}^{2}$ with $\mathbb{D}$, respectively).

In order to introduce the definition of the generalized Bloch space, we need the following notion.

A continuous increasing function $\omega\colon [0,\infty)\rightarrow[0,\infty)$ with $\omega(0)=0$ is called a {\it majorant} if $\omega(t)/t$ is non-increasing for $t>0$ (cf. \cite{dyak1997, dyak2004, Pa2004, Pa2007}).

 Given a subset $\Omega$ of $\mathbb{R}^{n}$, a function $f:\Omega\rightarrow\mathbb{R}^{n} $ is said to belong to the {\it Lipschitz space $\mathcal{L}_{\omega}(\Omega)$} if there is a positive constant $\mu_{01}$ such that for all $x,$ $y\in \Omega$,
$$|f(x)-f(y)|\leq \mu_{01}\omega(|x-y|).$$

First, we define the Bloch space of $\mathbb{B}^{n}$, denoted by $\mathcal{B}(\mathbb{B}^{n})$, as the space of functions $f$ in $C^{1}(\mathbb{B}^{n},\mathbb{R}^{n})$ such that $$||f||^{\mathcal{B}}<\infty,$$
where $$||f||^{\mathcal{B}}=\sup_{x\in\mathbb{B}^{n}}\big\{||D f(x)||(1-|x|^{2})\big\}$$
and $||D f(x)||$ denotes the matrix norm of the usual Jacobian matrix $D f(x)$ of $f$ at $x$ (See $\S$\ref{Matrix} below for the precise definition of $||D f(x)||$).
This is only a semi-norm. Obviously, $||f||^{\mathcal{B}}=0$ if and only if $f$ is constant.
The Bloch space $\mathcal{B}(\mathbb{B}^{n})$ becomes a Banach space with the following norm:
$$||f||_{\mathcal{B}}=|f(0)|+||f||^{\mathcal{B}}.$$

For an analytic function $f$ in $\mathbb{D}$, obviously, $||D f(z)||=|f'(z)|$.
Therefore, the classical {\it Bloch space $\mathcal{B}(\mathbb{D} )$} consisting of the analytic functions $f$ satisfying
 $$||f||_{\mathcal{B}}=|f(0)|+\sup_{x\in\mathbb{D}}\big\{|f'(z)|(1-|z|^{2})\big\}<\infty$$
is a subspace of $\mathcal{B}(\mathbb{B}^{n})$ (cf. \cite{anclpo, howa, zhu2005}).

Next, we introduce the notion: {\it Generalized Bloch spaces}.

\bdefe\label{def-1.3}
For $p\in (0,\infty]$, $\alpha>0$, $\beta\in \mathbb{R}$ and a majorant $\omega$, we use
$\mathcal{L}_{p,\omega} \mathcal{B}_{\alpha, \mathrm{a}}^{\beta}(\mathbb{B}^{n})$ to denote the {\it generalized Bloch space}, which consists of all functions $f\in C^{1}(\mathbb{B}^{n},\mathbb{R}^{n})$ with $$||f||_{\mathcal{L}_{p,\omega} \mathcal{B}_{\alpha, \mathrm{a}}^{\beta}(\mathbb{B}^{n})}<\infty,$$
where
$$||f||_{\mathcal{L}_{p,\omega} \mathcal{B}_{\alpha, \mathrm{a}}^{\beta}(\mathbb{B}^{n})}=|f(0)|+\sup_{x\in \mathbb{B}^{n}}\Big\{M_{p}(|x|,||Df||)\omega\big(\phi_{\alpha, \beta, \mathrm{a}}(x)\big)  \Big\},
$$
$$\phi_{\alpha, \beta, \mathrm{a}}(x)=(1-|x|)^{\alpha}\left( \log\frac{\mathrm{a}}{1-|x|}\right)^{\beta}$$
and $\mathrm{a}$ is a constant satisfying
\begin{enumerate}
  \item $\mathrm{a}>1$, if $\beta\leq 0$; and
  \item $\mathrm{a}\geq e^{\frac{\beta}{\alpha}}$, if $\beta>0$.
\end{enumerate}
\edefe

The space $\mathcal{L}_{p,\omega} \mathcal{B}_{\alpha, e}^{\beta}(\mathbb{D})$ was discussed in \cite{chensami2015} and the space $\mathcal{L}_{p,id} \mathcal{B}_{\alpha, e^{\beta/\alpha}}^{\beta}(\mathbb{D})$ of analytic functions was introduced in \cite{stev09}. Note that, when $\beta=0$, the space $\mathcal{L}_{p,\omega} \mathcal{B}_{\alpha, \mathrm{a}}^{0}(\mathbb{B}^{n})$ has nothing to do with the parameter $\mathrm{a}$.

Obviously, $\mathcal{L}_{\infty,id} \mathcal{B}_{1,\mathrm{a}}^{0}(\mathbb{B}^{n})=\mathcal{B}(\mathbb{B}^{n})$, where $id$ denotes the identity mapping. Further, we have the following.
\begin{enumerate}
\item
  The special case $\mathcal{L}_{\infty,\omega} \mathcal{B}_{\alpha, \mathrm{a}}^{0}(\mathbb{D})$ is called the {\it $\omega$-$\alpha$-Bloch space}
  (cf. \cite{fu2016, yoneda} and the related references therein).
\item
 The special case $\mathcal{L}_{\infty,\omega} \mathcal{B}_{1,\mathrm{a}}^{\beta}(\mathbb{D})$ is called the {\it logarithmic $\omega$-Bloch space}
 (cf. \cite{stev08, stev09} and the related references therein).
 \item
 The special case $\mathcal{L}_{\infty,id} \mathcal{B}_{\alpha,\mathrm{a}}^{0}(\mathbb{D})$ is called the {\it generalized $\alpha$-Bloch space}
 (cf. \cite{liwu, liste1, liste} and the related references therein).
 \item
 The special case $\mathcal{L}_{\infty,id} \mathcal{B}_{1,\mathrm{a}}^{\beta}(\mathbb{D})$ is called the {\it generalized logarithmic Bloch space}
 (cf. \cite{petr2013, yoneda} and the related references therein).
\end{enumerate}

 In \cite{dyak1997}, Dyakonov discussed the relationship between the Lipschitz space $\mathcal{L}_{\omega}(\mathbb{D})$ and the bounded mean oscillation on analytic functions in $\mathbb{D}$ (\cite[Theorem 1]{dyak1997}). In \cite{chensami2015} and \cite{chenrasi2014}, the authors extended \cite[Theorem 1]{dyak1997} to the case of complex-valued harmonic mappings (\cite[Theorem 4]{chensami2015} and \cite[Theorem 3]{chenrasi2014}). Recently, Chen and Rasila generalized \cite[Theorem 4]{chensami2015}, \cite[Theorem 3]{chenrasi2014} and \cite[Theorem 1]{dyak1997} to the setting of the solutions to the non-homogenous Yukawa PDE
$\Delta f=\lambda f$, where $\lambda:\mathbb{B}^{n}\rightarrow \mathbb{R}$ is a nonnegative continuous function with $\sup_{x\in\mathbb{B}^{n} }\{\lambda(x)\}<\infty$ (\cite[Theorem 1]{chenantti2016}).
 As the first aim of this paper, we consider the similar results of the above type for hyperbolic harmonic mappings. The following is our first result in this line.

 \begin{thm}\label{thm-1.1}
Suppose $\alpha\in [1,2)$, $\omega$ is a majorant and $\Delta_{h} u=0$. Then
$u\in \mathcal{L}_{\infty,\omega}\mathcal{B}_{\alpha,\mathrm{a}}^{0}(\mathbb{B}^{n})$
if and only if there is a positive constant $\mu_1$
 such that for all $r\in (0,1-|x|)$,
%\be\label{eq-1.2}
$$\frac{1}{|\mathbb{B}^{n}(x,r)|} \int_{\mathbb{B}^{n}(x,r)}|u(y)-u(x)|d\nu(y)\leq \frac{\mu_1}{\omega(r^{\alpha})} r,$$
where $|\mathbb{B}^{n}(x,r)|$ means the Lebesgue volume of the ball $\mathbb{B}^{n}(x,r)$ and $d\nu$ denotes the normalized Lebesgue volume measure in $\mathbb{B}^{n}$.
\end{thm}
\br
If we take $n=2$, $\alpha=1$ and replace $\omega(t)$ by $\frac{1}{\omega(\frac{1}{t})}$ in Theorem \ref{thm-1.1}, then
Theorem \ref{thm-1.1} coincides with \cite[Theorem 3]{chenrasi2014}.
\er

To state our next result, let us recall the following notion.

The {\it hyperbolic distance} between two points $x$ and $y$ in $\mathbb{B}^{n}$ is defined by $$\rho(x,y)=\inf_{\gamma\in\Gamma_{xy}(\mathbb{B}^{n})}\int_{\gamma}\frac{2}{1-|z|^{2}} ds(z) ,$$
where $ds$ is the length element on $\gamma$ and $\Gamma_{xy}(\mathbb{B}^{n})$ stands for the collection of all rectifiable curves in $\mathbb{B}^{n}$ joining $x$ and $y$ (cf. \cite{matti1988}).
See $\S$\ref{hyperbolic-metric} below for more properties of $\rho$.

In \cite{zhu2005}, Zhu characterized the holomorphic Bloch space in $\mathbb{C}^{n}$ in terms of the Bergman metric (\cite[Theorem 3.6 and Corollary 3.7]{zhu2005}).
Now, we establish the following necessary and sufficient condition for hyperbolic harmonic mappings to be in $\mathcal{B}(\mathbb{B}^{n})$ in terms of the hyperbolic metric.

\begin{thm}\label{thm-1.2}
Suppose $\Delta_{h} u=0$. Then $u\in \mathcal{B}(\mathbb{B}^{n})$ if and only if there exists a positive constant $\mu_{2}$ such that for all $x,y\in \mathbb{B}^{n}$,
\be\label{eq-1.3}|u(x)-u(y)|\leq \mu_{2}\rho(x,y).\ee
 \end{thm}

%%%%%%%%%%%%%%%%%%%%%%%%%%%%%%%%%%%%
\subsection{Integral means}
%%%%%%%%%%%%%%%%%%%%%%%%%%%%%%%%%%%%

First, we recall the following well known result on analytic functions due to Hardy and Littlewood.

\begin{Thm}\label{Thm-O} $($\cite[Theorems 2 and 3]{hard1931} {\rm and}  \cite[Theorem 46]{hard1932} {\rm or} \cite[Theorem 5.5]{dure}{\rm )}
Suppose $p\in (0,\infty]$, $\alpha\in (1,\infty)$ and $f$ is an analytic function in $\mathbb{D}$. Then the following statements are equivalent.
\ben
\item
$M_{p}(r,f')=O\left(\frac{1}{(1-r)^{\alpha}} \right)$ as $r\rightarrow 1 $;
\item
$M_{p}(r,f)=O\left(\frac{1}{(1-r)^{\alpha-1}} \right)$ as $r\rightarrow 1.$
\een
\end{Thm}

Obviously, the above result of Hardy and Littlewood provides a close relationship between the integral means of analytic functions and that of their derivatives.

As the second aim of this paper, we consider Theorem \Ref{Thm-O} in the setting of hyperbolic harmonic mappings. Our first result is the following analog of the implication from (1) to (2) in Theorem \Ref{Thm-O} for hyperbolic harmonic mappings.

\begin{thm}\label{thm-1.5}
Suppose $\Delta_{h}u=0$ and $u\in \mathcal{L}_{p,\omega}\mathcal{B}_{\alpha,\mathrm{a}}^{\beta}(\mathbb{B}^{n})$.
Then for $r\in [0,1)$ and $p\in [1,\infty]$,
  $$M_{p}(r,u)\leq |u(0)|+\frac{(\log\mathrm{a})^{\beta} ||u||_{\mathcal{L}_{p,\omega}\mathcal{B}_{\alpha,\mathrm{a}}^{\beta}(\mathbb{B}^{n}) } }{ \omega \big((\log\mathrm{a})^{\beta}\big) } \int_{0}^{1} \frac{ r}{ \phi_{\alpha, \beta, \mathrm{a}}( r t)}dt .$$
 \end{thm}

 \br
 By taking $n=2$, $\omega(t)=t$, $\alpha>1$ and $\beta=0$, we obtain that Theorem \ref{thm-1.5} is a generalization of the implication from $(1)$ to $(2)$ in Theorem \Ref{Thm-O} even in the case of harmonic mappings when $p\in [1,\infty]$.
\er

We also consider the converse of Theorem \ref{thm-1.5}, and we get the following analog of the implication from $(2)$ to $(1)$ in Theorem \Ref{Thm-O} for hyperbolic harmonic mappings.
\begin{thm}\label{thm-1.6}
Suppose $\Delta_{h} u=0$ and
$$ M_{p}(r,u) =O\left(\frac{1}{(1-r)^{\alpha}}\right)\;\;{\rm as}\;\; r\rightarrow 1,$$ where $p\in (0,\infty)$ and $\alpha>\frac{1}{p}$. Then
for $q\in (0,\infty]$,
$$u\in \mathcal{L}_{q,id} \mathcal{B}_{\alpha+1+\frac{n-1}{p}, \mathrm{a}}^{0}(\mathbb{B}^{n}).$$
 \end{thm}

%%%%%%%%%%%%%%%%%%%%%%%%%%%%%%%%%%%%
\subsection{Weak uniform boundedness property}
%%%%%%%%%%%%%%%%%%%%%%%%%%%%%%%%%%%%

Let $\Omega$ be a proper domain of $\mathbb{R}^{n}$.
For $x\in \Omega$, we use $d_{\Omega}(x)$ to denote the Euclidean distance from $x$ to the boundary $\partial \Omega$ of $\Omega$.
For $x,y\in \Omega$, let
$$r_{\Omega}(x,y)=\frac{|x-y|}{\min\{d_{\Omega}(x),d_{\Omega}(y)\}}\;\;\text{and}
\;\;k_{\Omega}(x,y)=\inf_{\gamma\in\Gamma_{xy}(\Omega)}\int_{\gamma}\frac{1}{d_{\Omega}(z)}ds(z)$$
(cf. \cite{matel, matti1988}).

We say that $f:\Omega\rightarrow f(\Omega)\subset \mathbb{R}^{n}$ satisfies the {\it weak uniform boundedness property} in $\Omega$ (with respect to $r_{\Omega}$) if there is a constant $\mu_{02}>0$ such that for all $x,$ $y\in \Omega$,
\be\label{eq-1.4}r_{\Omega}(x,y)\leq \frac{1}{2}\;\; {\rm implies}\;\; r_{f(\Omega)}\big(f(x),f(y)\big)\leq \mu_{02}.\ee

\br
 The above definition of the weak uniform boundedness property is equivalent to the following one:
 $f:$ $\Omega\rightarrow f(\Omega)$ is said to satisfy the {\it weak uniform boundedness property} in $\Omega$ (with respect to $r_{\Omega}$) if for any $\mu_{03}\in (0,1)$, there is a constant $\mu_{04}>0$ such that for all $x,$ $y\in \Omega$,
$$r_{\Omega}(x,y)\leq \mu_{03}\;\; {\rm implies}\;\; r_{f(\Omega)}\big(f(x),f(y)\big)\leq \mu_{04}.$$
\er

In \cite[Theorem 2.8]{matel}, Mateljevi\'{c} and Vuorinen proved that a harmonic mapping $f$ satisfies the weak uniform boundedness property in $G\subset\mathbb{R}^{n}$ if and only if there exists a constant $\mu_{05}$ such that for all $x,y\in G$,
$$k_{f(G)}(f(x),f(y))\leq \mu_{05}k_{G}(x,y).$$
Recently, Chen and Rasila generalized this result to the case of the solutions to Yukawa equation in $\mathbb{B}^{n}$ (\cite[Theorem 2]{chenantti2016}).
As the last aim of this paper, we consider the weak uniform boundedness property of  hyperbolic harmonic mappings. Our result reads as follows.

\begin{thm}\label{thm-1.8}
Suppose $\Delta_{h} u=0$. Then $u$ satisfies the weak uniform boundedness property in $\mathbb{B}^{n}$ if and only if there exists a positive constant $\mu_3$ such that for all $x,y\in \mathbb{B}^{n}$,
$$k_{u(\mathbb{B}^{n})}\big(u(x),u(y)\big)\leq \mu_3k_{\mathbb{B}^{n}}(x,y). $$ %where $\mu_4=\mu_4(n,\mu_{02})$. Here $\mu_{02}$ is the constant from \eqref{eq-1.4}.
\end{thm}

This paper is organized as follows. In Section \ref{sec-2}, some necessary terminology and notation will be introduced.
In Section \ref{sec-3}, we shall prove Theorem \ref{thm-1.1}, and in Section \ref{sec-4}, we shall show Theorem \ref{thm-1.2}.
The proofs of Theorems \ref{thm-1.5} and \ref{thm-1.6} will be presented in Section \ref{sec-5}. Section \ref{sec-7} will be devoted to the proof of Theorem \ref{thm-1.8}.

%%%%%%%%%%%%%%%%%%%%%%%%%%%%%%%%%%%%%%%%%%%%%%%%
%%%%%%%%%%%%%%%%%%%%%%%%%%%%%%%%%%%%%%%%%%%%%%%%
\section{Preliminaries}\label{sec-2}
%%%%%%%%%%%%%%%%%%%%%%%%%%%%%%%%%%%%%%%%%%%%%%%%
%%%%%%%%%%%%%%%%%%%%%%%%%%%%%%%%%%%%%%%%%%%%%%%%

In this section, we recall some necessary terminology and notation.

%%%%%%%%%%%%%%%%%%%%%%%%%%%%%%%%%%%%%%%%%%%%%%%%
\subsection{Matrix notations}\label{Matrix}
%%%%%%%%%%%%%%%%%%%%%%%%%%%%%%%%%%%%%%%%%%%%%%%%

For a natural number $n$, let
$$A=\big(a_{ij}\big)_{n\times n}\in \mathbb{R}^{n\times n}. $$

For $A\in \mathbb{R}^{n\times n},$
denote by $\|A\|$ the matrix norm $$\|A\|=\sup\{|Ax|:\; x\in \mathbb{R}^{n}, |x|=1\}.$$% and by $l(A)$ the matrix function $$l(A)=\inf\{|Ax|:\; x\in \mathbb{R}^{n}, |x|=1\}.$$

For a domain $\Omega\subset\mathbb{R}^{n}$, let $f=(f_{1},\ldots,f_{n})\colon \Omega\rightarrow \IR^n$ be a function that has all partial derivatives at $x=(x_{1},\ldots,x_{n})$ in $\Omega$. Then $Df(x)$ denotes the usual Jacobian matrix
$$
Df=\left(
    \begin{array}{cccc}
      \frac{\partial f_{1}(x)}{\partial x_{1}}    & \cdots & \frac{\partial f_{1}(x)}{\partial x_{n}} \\
     \vdots & \ddots  & \vdots \\
    \frac{\partial f_{n}(x)}{\partial x_{1}}  & \cdots & \frac{\partial f_{n}(x)}{\partial x_{n}}\\
    \end{array}
  \right) = \big(\nabla f_1(x) \cdots \nabla f_n(x)\big)^T
$$
at $x$, where $T$ is the transpose and the gradients $\nabla f_j(x)$ are understood as column vectors (cf. \cite{chrw}).
For two column vectors $x$, $y\in \mathbb{R}^{n}$, we use $\langle x,y\rangle$ to denote the inner product of $x$ and $y$.
For $j\in\{1,\cdots,n\}$, it follows from $$||Df(x)||^{2}=\sup_{\xi\in \mathbb{S}^{n-1}}\Big\{\sum_{j=1}^{n}\big(\langle\nabla f_{j}(x),\xi\rangle\big)^{2}\Big\}$$ that
\be\label{eq-2.8}||Df(x)||^{2}\leq\sum_{j=1}^{n}|\nabla f_{j}(x)|^{2}=\sum_{i=1}^{n}\left|\frac{\partial }{\partial x_{i}} f(x)\right|^{2}\;\;\text{and}\;\;||Df(x)||\geq|\nabla f_{j}(x)|.\ee

%%%%%%%%%%%%%%%%%%%%%%%%%%%%%%%%%%%%%%%%%%%%%%%%
 \subsection{Hyperbolic metric}\label{hyperbolic-metric}
%%%%%%%%%%%%%%%%%%%%%%%%%%%%%%%%%%%%%%%%%%%%%%%%

 For any $w\in \mathbb{B}^{n}$, let
$$\varphi_{w}(x)=\frac{|x-w|^{2}w-(1-|w|^{2})(x-w)}{[x,w]^{2}}$$
 in $\mathbb{B}^{n}$, where $[x,w]=\left||x|w -\frac{x}{|x|}\right|$.
 Then $\varphi_{w}$ is a M\"{o}bius transformation from $\mathbb{B}^{n}$ onto $\mathbb{B}^{n}$
 with $\varphi_{w}(w)=0$, $ \varphi_{w}(0)=w$ and $\varphi_{w}\big(\varphi_{w}(x)\big)=x$.

We denote by $\mathcal{M}(\mathbb{B}^{n})$ the set of all M\"{o}bius transformations in $\mathbb{B}^{n}$. It is well known that
 if $\varphi\in \mathcal{M}(\mathbb{B}^{n})$, then there exist $w\in \mathbb{B}^{n}$ and
 an orthogonal transformation $A$ such that $$\varphi(x)=A\varphi_{w}(x)$$ (cf. \cite[Theorem 2.1]{sto1999} or \cite{sto2016}).

For more information on M\"{o}bius transformations in $\mathbb{B}^{n}$, see e.g. \cite{ahlf, beardon, matti1988}.

In terms of $ \varphi_{w}$, the {\it hyperbolic metric} $\rho$ in $\mathbb{B}^{n}$ can be given by
\be\label{eq-2.1}\rho(x,w)=\log \frac{1+|\varphi_{w}(x)|}{1-|\varphi_{w}(x)|}\ee
for $ x,$ $w\in \mathbb{B}^{n}$.
By \cite[Equation (2.6)]{sto1999} or \cite[Equation (2.11)]{chrw}, we see that for $x,w\in \mathbb{B}^{n}$,
\be\label{eq-2.10}|\varphi_{w}(x)|=|\varphi_{x}(w)|=\frac{|x-w|}{[x,w]}.\ee
 Therefore, $\rho(x,w)=\rho(w,x)$.
In particular,
$$\rho(0,x)=\log\frac{1+|x|}{1-|x|}  \;\;\text{and}\;\;  \rho(x,y)=\rho\big(\varphi(x),\varphi(y)\big)$$
for all $\varphi\in \mathcal{M}(\mathbb{B}^{n})$ (cf. \cite[Chapter 1]{matti1988}).

 For any $w\in \mathbb{B}^{n}$ and $0<r<1$, we define the {\it pseudo-hyperbolic ball} with center $w$ and radius $r$ as
 \be\label{eq-2.02}E(w,r)= \{x\in \mathbb{B}^{n}:| \varphi_{w}(x)|<r\}.\ee
Clearly, $E(w,r)= \varphi_{w}\big(\mathbb{B}^{n}(0,r)\big)$ (cf. \cite{sto2012}).

It is well known that $E(w,r)$ is also a Euclidean ball with centre $c_{w}$ and radius $r_{w}$ which are as follows:
  \be\label{eq-2.2}c_{w}=\frac{1-r^2}{1-|w|^{2}r^{2}}w\;\;\;\text{and}\;\;\;r_{w}=\frac{1-|w|^2}{1-|w|^{2}r^{2}}r.\ee

By \cite[Lemma 2.1]{reb2005}, we obtain that for any $\delta\in(0,1)$ and $y\in E(x,\delta)$,
$$1-|y|^2\leq \left||x|y-\frac{x}{|x|}\right|(1+|y|)\leq\frac{2(1+\delta)}{1-\delta}(1-|x|^2).$$
By \eqref{eq-2.10} and \eqref{eq-2.02}, we see that $y\in E(x,\delta)$ if and only if $x\in E(y,\delta)$.
Therefore,
\be\label{eq-2.3}1-|x|^2\leq \frac{2(1+\delta)}{1-\delta}(1-|y|^2).
\ee

%%%%%%%%%%%%%%%%%%%%%%%%%%%%%%%%%%%%%%%%%%%%%%%%
\subsection{Hyperbolic harmonic mappings}\label{HyperP}
%%%%%%%%%%%%%%%%%%%%%%%%%%%%%%%%%%%%%%%%%%%%%%%%

For all $\varphi \in \mathcal{M}(\mathbb{B}^{n})$ and $f\in C^{2}(\mathbb{B}^{n},\mathbb{R}^{n})$, we have the following M\"obius invariance property (cf. \cite[Section 2]{sto2012}):
\be\label{eq-2.4}\Delta_{h}( f\circ \varphi )(x)=\Delta_{h} f\big(\varphi(x)\big) .\ee
Obviously, \eqref{eq-01.01} implies that
$$ \Delta_{h} f (x)= \Delta_{h}f\big(\varphi_{x}(0)\big)=\Delta_{h}(f\circ \varphi_{x})(0)= \Delta(f\circ \varphi_{x})(0).$$

It is well known that if $\psi\in C(\mathbb{S}^{n-1},\mathbb{R}^{n})$, then the Dirichlet problem
$$\begin{cases}
\displaystyle \;\Delta_{h} f=0\;\;
\text{in }\;\;\mathbb{B}^n,\\
\displaystyle \;f=\psi\;\;\;\;\;\; \text{on }\;\;\mathbb{S}^{n-1}
\end{cases}
$$
has a unique solution in $C(\overline{\mathbb{B}}^{n})$ and can be represented by
$$f(x)=P_{h}[\psi](x)=\int_{\mathbb{S}^{n-1}}P_{h}(x,\xi)\psi(\xi)\,d\sigma(\xi)$$
(cf. \cite{chenshao2015} or \cite{sto1999}), where
$$ P_{h}(x,\xi)=\left(\frac{1-|x|^2}{|x-\xi|^2}\right)^{n-1}.$$

For $f\in C^{1}(\mathbb{B}^{n},\mathbb{R})$, the gradient $\nabla^{h}f$ with respect to the hyperbolic metric is given by
$$\nabla^{h}f(x)=(1-|x|^{2})\nabla f(x)$$
(cf. \cite{pav1991, reka, sto1999}).
In particular,
\be\label{eq-2.5}|\nabla^{h}f(x)|=(1-|x|^{2})|\nabla f(x)|.\ee
 Furthermore, for all $\varphi\in \mathcal{M}(\mathbb{B}^{n}) $, by \cite[Theorem 3.2]{sto1999} or \cite[Equation (2.14)]{sto2012}, we have
\be\label{eq-2.6} |\nabla^{h}(f\circ \varphi)(x)|=\big|\nabla^{h} f\big( \varphi(x)\big)\big| .\ee

From \cite[Lemma 3.3]{sto2012} or \cite[Theorem 2.1]{pav1991}, we can easily obtain the following useful result.

\begin{Lem}\label{lem-B}
Suppose $p\in (0, \infty)$ and $\delta\in (0, \frac{1}{2})$. Then there exists a constant $\mu_{06}$ such that for any $f\in C^{2}(\mathbb{B}^{n},\mathbb{R}^{n})$ with $\Delta_{h} f=0$,
$$ |\nabla^{h} f(x)|^{p}\leq \mu_{06}\delta^{-n}\int_{E(x,\delta)} | f(y)|^{p}\;d\tau(y)$$
 in $\mathbb{B}^{n}$, where $\mu_{06}=\mu_{06}(p,\delta)$ (which means that the constant $\mu_{06}$ depends only on the given parameters $p$ and $\delta$) and $d\tau$ denotes the {\it M\"{o}bius invariant measure} in $\mathbb{B}^{n}$, which is given by
%\be\label{eq-2.7}
$$d\tau(x)=\frac{d\nu(x)}{(1-|x|^2)^n}.$$
\end{Lem}

%%%%%%%%%%%%%%%%%%%%%%%%%%%%%%%%%%%%
%%%%%%%%%%%%%%%%%%%%%%%%%%%%%%%%%%%%
\section{Generalized Bloch spaces and bounded mean oscillation}\label{sec-3}
%%%%%%%%%%%%%%%%%%%%%%%%%%%%%%%%%%%%
%%%%%%%%%%%%%%%%%%%%%%%%%%%%%%%%%%%%

The purpose of this section is to prove Theorem \ref{thm-1.1}. Before the proof, we need the following lemma.
\begin{lem}\label{lem-3.1}
Suppose $\Delta_{h} u=0$.
Then there is a constant $\mu_{07}$ such that
for any $x\in \mathbb{B}^{n}$,
$$||Du(x)||\leq \frac{5^{n}\sqrt{n}\mu_{07}}{2^{n}(1-|x|)^{n+1}}\int_{\mathbb{B}^{n}\left(x,\frac{1}{4}(1-|x|)\right)}|u(y)-u(x)|\;d\nu(y),$$
where $\mu_{07}$ is a constant depending only on $\mu_{06}(1,\frac{1}{9})$ and $\mu_{06}$ is the constant from Lemma \Ref{lem-B}.
\end{lem}
\bpf
For $x\in \mathbb{B}^{n}$, by \eqref{eq-2.8}, we know that
\be\label{nd-3}
\|Du(x)\|
\leq\left(\sum_{j=1}^{n}\left|\nabla u_{j}(x)\right|^{2}\right)^{\frac{1}{2}}.\ee
Hence, to prove this lemma, it suffices to estimate the quantity $\left|\nabla u_{j}(x)\right|$.
Since for a fixed $x\in \mathbb{B}^n$, $$\Delta_{h}\big(u(w)-u(x)\big)=0$$ in $\mathbb{B}^n$, by taking $p=1$ and $\delta=\frac{1}{9}$, we see from Lemma \Ref{lem-B} and \eqref{eq-2.5} that there is a constant $\mu_{08}$ such that
for any $w\in \mathbb{B}^n$ and each $j\in\{1,\cdots,n\}$,
\beq\label{nd-31}
%\begin{eqnarray*}
(1-|w|^{2})|\nabla u_{j}(w)|
&\leq& \mu_{08}\int_{E(w,\frac{1}{9})} |u_{j}(y)-u_{j}(x)|d\tau(y)\\  \nonumber
&\leq& \mu_{08}\int_{E(w,\frac{1}{9})} |u(y)-u(x)|d\tau(y).
\eeq
Therefore, \eqref{nd-3} guarantees the following:
\beq\label{fri-11}
%\begin{eqnarray*}
(1-|w|^{2})||Du(w)||
&\leq&(1-|w|^{2})\left(\sum_{j=1}^{n}|\nabla u_{j}(w) |^{2}\right)^{\frac{1}{2}}\\  \nonumber
&\leq& \sqrt{n}\mu_{08}\int_{E(w,\frac{1}{9})} |u(y)-u(x)|d\tau(y)\\  \nonumber
&=&\sqrt{n} \mu_{08}\int_{E(w,\frac{1}{9})} \frac{|u(y)-u(x)|}{(1-|y|^{2})^{n}}d\nu(y).
\eeq
Moreover, by \eqref{eq-2.2}, we know that
\be\label{fri-21}E(x,\frac{1}{9})= \mathbb{B}^{n}\left(\frac{80}{81-|x|^{2}}x,\frac{9(1-|x|^{2})}{81-|x|^{2}}\right)\subset  \mathbb{B}^{n}\left( x,\frac{1}{4}(1-|x| )\right).\ee
Further, for any $y\in E(x,\frac{1}{9})$, \eqref{eq-2.3} implies
\be\label{fri-22}
\frac{1}{1-|y|^{2}}\leq\frac{5}{2(1-|x|)}.
\ee
By letting  $w=x$,
it follows from \eqref{fri-11}$\sim$\eqref{fri-22} that
$$
(1-|x|^{2}) ||Du(x)|| \leq \frac{5^{n} \sqrt{n} \mu_{08}}{2^{n}(1-|x|)^{n}}\int_{ \mathbb{B}^{n} \left(x,\frac{1}{4}(1-|x|)\right)}  |u(y)-u(x)| d\nu(y).
$$
By taking $\mu_{07}= \mu_{08}$, we know that the lemma is proved.
\epf

%%%%%%%%%%%%%%%%%%%%%%%%%%%%%%%%%%%%
\subsection*{Proof of Theorem \ref{thm-1.1}}
%%%%%%%%%%%%%%%%%%%%%%%%%%%%%%%%%%%%

First, we show the ``if" part in the theorem.
For any $x\in \mathbb{B}^{n}$, we see that
$$ \frac{1}{2^{n}(1-|x|)^{n+1} }=\frac{1}{8^{n}(1-|x|) \big(\frac{1}{4}(1-|x| )\big)^{n} }=\frac{ |\mathbb{B}^{n}|}{8^{n}(1-|x|) \left|\mathbb{B}^{n}\left(x,\frac{1}{4}(1-|x|)\right)\right| }. $$
Then by Lemma \ref{lem-3.1}, we know that there is a constant $\mu_{07}$ such that for any $x\in \mathbb{B}^{n}$,
$$||Du(x)||\leq \frac{5^{n}\sqrt{n}\mu_{07}|\mathbb{B}^{n}|}{8^{n}(1-|x|) \left|\mathbb{B}^{n}\left(x,\frac{1}{4}(1-|x|)\right)\right|}\int_{\mathbb{B}^{n}\left(x,\frac{1}{4}(1-|x|)\right)}|u(y)-u(x)|\;d\nu(y).$$
Then the assumption in the theorem implies
$$||Du(x)|| \leq \frac{ 5^{n}\sqrt{n} \mu_1\mu_{07}}{2^{3n+2}\omega\big(\frac{(1-|x|)^{\alpha} }{4^{\alpha}}\big)}|\mathbb{B}^{n}|.$$

Moreover, it follows from $\alpha\in[1,2) $ and \cite[Lemma 2.2(2)]{chenshao2015} that
$$ \omega\left(4^{\alpha}\cdot \frac{ (1-|x|)^{\alpha}}{4^{\alpha}}\right) \leq 4^{\alpha}\omega\left(\frac{ (1-|x|)^{\alpha}}{4^{\alpha}}\right),$$
and thus, for all $x\in \mathbb{B}^{n}$, $$|| Du(x) || \leq  \frac{ 5^{n}\sqrt{n} \mu_1\mu_{07}}{2^{3n+2-2\alpha}\omega\big( (1-|x|)^{\alpha} \big)}|\mathbb{B}^{n}|,$$
which means that $u\in\mathcal{L}_{\infty,\omega}\mathcal{B}_{\alpha,\mathrm{a}}^{0}(\mathbb{B}^{n})$.\medskip

Next, we prove the ``only if" part. Let $x\in \IB^n$.
For any $r\in (0, 1-|x|)$ and $y\in \mathbb{B}^{n}(x,r)$, obviously, we have
\be\label{eq-3.4}|y-x|<r< 1-|x|\;\;{\rm and}\;\; t|x-y|<1-|x|,\ee
 where $t\in[0,1]$. For the proof, we need an upper bound on the quantity $|u(y)-u(x)|$.  For this, we let $\gamma_{[x,y]}$ denote the segment between $x$ and $y$ with the parametrization $\gamma(t)=(1-t)x+ty$, where $t\in[0,1]$. By the well-known gradient theorem (see, e.g.
\cite[Theorem 6.24]{wrudin}), we have that for each $j\in \{1,\ldots, n\}$,
$$\int_{\gamma_{[x,y]}} \langle \nabla u_j(\gamma), d\gamma\rangle
= \int_{0}^{1} \big\langle\nabla u_j \big(\gamma(t)\big), \gamma'(t)\big\rangle \,dt=\int_{0}^{1}   \frac{d}{d t} \big(u_j \circ \gamma (t)\big)\,dt = u_j(y)-u_j(x).$$
Note that
\[
Du\big(\gamma(t)\big)\times \gamma'(t) = \left(
    \begin{array}{cccc} \big\langle\nabla u_1 \big(\gamma(t)\big), \gamma'(t)\big\rangle\\
    \vdots \\
    \big\langle\nabla u_n \big(\gamma(t)\big), \gamma'(t)\big\rangle
    \end{array}\right),
\]
where $A\times B$ denotes the product of two matrices $A$ and $B$.
Hence
$$ u(y)-u(x)=\int_{0}^{1}  Du\big(\gamma(t)\big)\times \gamma'(t)\,dt,$$
and therefore
\be\label{eq-3.7}|u(y)-u(x)|=\left|\int_{0}^{1}  Du\big(\gamma(t)\big)\times \gamma'(t)\,dt\right|
\leq \int_{\gamma_{[x,y]}}   \|D u(\gamma) \|\cdot|d\gamma|.\ee

Moreover, it follows from the assumption $u\in \mathcal{L}_{\infty,\omega}\mathcal{B}_{\alpha,\mathrm{a}}^{0}(\mathbb{B}^{n})$ that for $x\in \mathbb{B}^{n}$, $$||Du(x)||\leq \frac{1}{\omega \big((1-|x|)^{\alpha} \big)}||u||_{\mathcal{L}_{\infty,\omega}\mathcal{B}_{\alpha,\mathrm{a}}^{0}(\mathbb{B}^{n})}.$$
We infer from the easy fact $|\gamma'(t)|=|x-y|$ that
%\beq\label{eq-3.5}
\begin{eqnarray*}
|u(y)-u(x)| &\leq & ||u||_{\mathcal{L}_{\infty,\omega}\mathcal{B}_{\alpha,\mathrm{a}}^{0}(\mathbb{B}^{n})} \int_{0}^{1} \frac{|x-y|}{\omega \big((1-|(1-t)x+ty|)^{\alpha} \big)}dt\\ \nonumber
&\leq&  ||u||_{\mathcal{L}_{\infty,\omega}\mathcal{B}_{\alpha,\mathrm{a}}^{0}(\mathbb{B}^{n})} \int_{0}^{1} \frac{|x-y|}{\omega \big((1-|x|-t|x-y|)^{\alpha} \big) }dt\\ \nonumber
&=&  ||u||_{\mathcal{L}_{\infty,\omega}\mathcal{B}_{\alpha,\mathrm{a}}^{0}(\mathbb{B}^{n})} \int_{0}^{|x-y|} \frac{1}{\omega \big((1-|x|-t)^{\alpha} \big) }dt,
\end{eqnarray*}
 since the inequalities \eqref{eq-3.4} and the assumption $\alpha\in [1,2)$ guarantee that
 \begin{eqnarray*}
 \omega \big((1-|(1-t)x+ty|)^{\alpha} \big)
% &=&\omega \big((1-|(x+t(y-x)|)^{\alpha} \big)\\
&\geq& \omega \big((1-|x|-t|y-x|)^{\alpha} \big).
\end{eqnarray*}
This is our needed upper bound on $|u(y)-u(x)|$.\medskip

Now, we are ready to finish the proof. Let $x-y=\eta\in \mathbb{B}^{n}$. A similar argument as in the proof of \cite[Theorem 1]{chenantti2016} implies that
\begin{eqnarray*}
\;\;\;\;\;\;&\;\;& \frac{1}{|\mathbb{B}^{n}(x,r)|} \int_{\mathbb{B}^{n}(x,r)}|u(y)-u(x)|d\nu(y)\\
&\leq&\frac{||u||_{\mathcal{L}_{\infty,\omega}\mathcal{B}_{\alpha,\mathrm{a}}^{0}(\mathbb{B}^{n})} }{|\mathbb{B}^{n}(x,r)|}
\int_{\mathbb{B}^{n}(x,r)}\left(\int_{0}^{|x-y|} \frac{1}{\omega \big((1-|x|-t)^{\alpha} \big) }dt \right)d\nu(y )\\
&=& \frac{||u||_{\mathcal{L}_{\infty,\omega}\mathcal{B}_{\alpha,\mathrm{a}}^{0}(\mathbb{B}^{n})} }{|\mathbb{B}^{n}(0,r)|}
\int_{\mathbb{B}^{n}(0,r)}\left(\int_{0}^{|\eta|} \frac{1}{\omega \big((1-|x|-t)^{\alpha} \big) }dt \right)d\nu(\eta) \\
&=& \frac{n||u||_{\mathcal{L}_{\infty,\omega}\mathcal{B}_{\alpha,\mathrm{a}}^{0}(\mathbb{B}^{n})} }{r^{n}} \int_{0}^{r}\rho^{n-1}\left\{ \int_{0}^{\rho}   \frac{1}{\omega \big((1-|x|-t)^{\alpha} \big) }dt \right\}d\rho\\
&\leq& \frac{n||u||_{\mathcal{L}_{\infty,\omega}\mathcal{B}_{\alpha,\mathrm{a}}^{0}(\mathbb{B}^{n})} }{r^{n}} \int_{0}^{r} \left( \int_{t}^{r}\rho^{n-1} d\rho \right)  \frac{1}{\omega \big((r-t)^{\alpha} \big) }dt  \\
&=& \frac{||u||_{\mathcal{L}_{\infty,\omega}\mathcal{B}_{\alpha,\mathrm{a}}^{0}(\mathbb{B}^{n})} }{r^{n}}
\int_{0}^{r} \frac{(r-t)(r^{n-1}+r^{n-2}t+\cdots+t^{n-1})}{\omega \big((r-t)^{\alpha} \big) } dt  \\
&\leq& \frac{n||u||_{\mathcal{L}_{\infty,\omega}\mathcal{B}_{\alpha,\mathrm{a}}^{0}(\mathbb{B}^{n})} }{r}
\int_{0}^{r} \frac{(r-t)}{\omega \big((r-t)^{\alpha} \big) } dt  \\
&\leq& \frac{n||u||_{\mathcal{L}_{\infty,\omega}\mathcal{B}_{\alpha,\mathrm{a}}^{0}(\mathbb{B}^{n})} }{r}
\int_{0}^{r} \frac{(r-t)^{\alpha}}{\omega \big((r-t)^{\alpha} \big) } (r-t)^{1-\alpha}dt  \\
&\leq& \frac{nr^{\alpha-1}||u||_{\mathcal{L}_{\infty,\omega}\mathcal{B}_{\alpha,\mathrm{a}}^{0}(\mathbb{B}^{n})} }{\omega(r^{\alpha})}
\int_{0}^{r} (r-t)^{1-\alpha}dt  \\
&\leq& \frac{n||u||_{\mathcal{L}_{\infty,\omega}\mathcal{B}_{\alpha,\mathrm{a}}^{0}(\mathbb{B}^{n})} }{(2-\alpha)}\frac{r}{\omega(r^{\alpha})},
\end{eqnarray*} which is what we need.
\qed

%%%%%%%%%%%%%%%%%%%%%%%%%%%%%%%%%%%%
%%%%%%%%%%%%%%%%%%%%%%%%%%%%%%%%%%%%
\section{Bloch space and hyperbolic metric}\label{sec-4}
%%%%%%%%%%%%%%%%%%%%%%%%%%%%%%%%%%%%
%%%%%%%%%%%%%%%%%%%%%%%%%%%%%%%%%%%%

The aim of this section is to prove Theorem \ref{thm-1.2}. We start this section with a lemma.
\begin{lem}\label{lem-4.1}
Suppose $\Delta_{h} u=0$ and $u\in \mathcal{B}(\mathbb{B}^{n})$. Then for any $z\in \mathbb{B}^{n}$,
%\be\label{eq-4.1}
$$|u(z)-u(0)|\leq\frac{||u||^{\mathcal{B}}}{2} \log\frac{1+|z|}{1-|z|}.$$
\end{lem}

\bpf For any $z\in \mathbb{B}^{n}$, we let $\gamma_{[0,z]}$ denote the segment between $0$ and $z$ with the parametrization $\gamma(t)=tz$, where $t\in[0,1]$.
A similar argument as in \eqref{eq-3.7} leads to
 $$|u(z)-u(0)|\leq \int_{\gamma_{[0,z]}}   \|D u(\gamma) \|\cdot|d\gamma|=|z|\int_{0}^{1}  \|D u(tz) \|\; dt.$$
 Then the assumption $$ ||u||^{\mathcal{B}} =\sup_{x\in \mathbb{B}^{n}}(1-|x|^{2})||D u(x)||<\infty$$ implies that
$$|u(z)-u(0)|\leq ||u||^{\mathcal{B}} \int_{0}^{1} \frac{|z|}{1-|tz|^{2}} dt
=\frac{ ||u||^{\mathcal{B}}} {2} \log\frac{1+|z|}{1-|z|},$$
as needed.
\epf

%%%%%%%%%%%%%%%%%%%%%%%%%%%%%%%%%%%%
\subsection*{Proof of Theorem \ref{thm-1.2}}
%%%%%%%%%%%%%%%%%%%%%%%%%%%%%%%%%%%%

First, we show the ``if" part in the theorem.
Since $$||u||_{\mathcal{B}}= |u(0)|+\sup_{w\in \mathbb{B}^{n}} \big\{(1-|w|^{2})||D u(w)||\big\},$$
we know that, to prove this part, we need to estimate the quantity $(1-|w|^{2})||D u(w)||$.
Since for a fixed $x\in \mathbb{B}^n$,
$$\Delta_h\big(u(w)-u(x)\big)=0$$ in $\mathbb{B}^n$, by taking $p=1$ and $\delta=\frac{1}{9}$, we know from Lemma \Ref{lem-B}, \eqref{nd-3} and \eqref{nd-31} that there is a constant $\mu_{08}$ such that
for any $w\in \mathbb{B}^n$ and each $j\in\{1,\cdots,n\}$,
\begin{eqnarray*}
(1-|w|^{2})^{2}||Du(w)||^{2}
&\leq& (1-|w|^{2})^{2}\sum_{j=1}^{n}|\nabla u_{j}(w) |^{2}\\\nonumber
&\leq& \mu_{08}^{2}\sum_{j=1}^{n}\left(\int_{E(w,\frac{1}{9})} |u_{j}(y)-u_{j}(x)|d\tau(y)\right)^{2}.
\end{eqnarray*}
Further, H\"{o}lder inequality leads to
%\beq\label{fri-1}
\begin{eqnarray*}
&\;\;&(1-|w|^{2})^{2}||Du(w)||^{2}\\\nonumber
&\leq&\mu_{08}^{2}\left(\int_{E(w,\frac{1}{9})}d\tau(y)\right)\cdot\left(\int_{E(w,\frac{1}{9})} \sum_{j=1}^{n}|u_{j}(y)-u_{j}(x)|^{2}d\tau(y)\right) \\\nonumber
&=& \mu_{08}^{2}\tau\big(E(w,\frac{1}{9})\big)\int_{E(w,\frac{1}{9})} |u(y)-u(x)|^{2}d\tau(y)\\\nonumber
&\leq & \mu_{08}^{2} \left(\tau\big(E(w,\frac{1}{9})\big) \right)^{2}\sup_{y\in E(w,\frac{1}{9})}\big\{|u(y)-u(x)|^{2}\big\}.
\end{eqnarray*}
%\eeq
Moreover, by \cite[Equation (2.7)]{sto2012}, we know
$$\tau\Big( E(w,\frac{1}{9}) \Big)=n\int_{0}^{\frac{1}{9}}t^{n-1}(1-t^2)^{-n}dt\leq \frac{9^{n}}{80^{n}}.$$
By letting  $w=x$,
we see that
\begin{eqnarray*}
(1-|x|^{2})^{2}||Du(x)||^{2}&\leq & \mu^{2}_{2} \mu^{2}_{08}\frac{81^{n}}{80^{2n}} \sup_{y\in E(x,\frac{1}{9})}\big\{\rho^{2}(y,x)\big\}
\;\;\;\;\;\;\;\;\;\;\;\;\;\;\;\;\;\;\;\;\;\;\;\;\text{(by \eqref{eq-1.3})}\\
&=& \mu^{2}_{2} \mu^{2}_{08}\frac{81^{n}}{80^{2n}} \sup_{y\in E(x,\frac{1}{9})} \left\{ \log^{2}\left(\frac{1+| \varphi_{x}(y)|}{1-|\varphi_{x}(y)|}\right)\right\}
\;\;\;\;\;\text{(by \eqref{eq-2.1})}\\
&\leq&   \mu^{2}_{2} \mu^{2}_{08}\frac{81^{n}}{80^{2n}}\log^{2}\frac{10}{8}. \;\;\;\;\;\;\;\;\;\;\;\;\;\;\;\;\;\;\;\;\;\;\;\;\;\;\;\;\;\;\;\;\;\;\;\;\;\;\;\;\text{(by \eqref{eq-2.02})}
\end{eqnarray*}
Then the arbitrariness of $x$ in $\mathbb{B}^n$ ensures the following: $$||u||_{\mathcal{B}}= |u(0)|+\sup_{x\in \mathbb{B}^{n}} \big\{(1-|x|^{2})||D u(x)||\big\}
\leq |u(0)|+  \mu_{2}\mu_{08}\frac{9^{n}}{80^{n}} \log\frac{10}{8}  <\infty. $$\medskip

Next, we prove the ``only if" part. We shall show this part by applying Lemma \ref{lem-4.1} to the mapping $u\circ \varphi_{y}$, where $y\in \mathbb{B}^{n}$. Hence we have to verify that $u\circ \varphi_{y}$ satisfies the conditions in Lemma \ref{lem-4.1}. Obviously, the hyperbolic harmonicity of $u\circ \varphi_{y}$ easily follows from \eqref{eq-2.4}. It remains to check that $u\circ \varphi_{y}$ satisfies the second assumption in Lemma \ref{lem-4.1}, which is stated in the following claim.
\bcl\label{nd-4}
For any $y\in \mathbb{B}^{n}$, $u\circ \varphi_{y}\in \mathcal{B}(\mathbb{B}^{n})$.
\ecl

Since $u\circ \varphi_{y}\in \mathcal{B}(\mathbb{B}^{n})$ if and only if
$$||u\circ \varphi_{y}||^{\mathcal{B}}=\sup_{x\in\mathbb{B}^{n}}\{||D (u\circ \varphi_{y})(x)||(1-|x|^{2})\}<\infty,$$ clearly, to check this claim, it needs to estimate the quantity
$$||D (u\circ \varphi_{y})(x)||(1-|x|^{2}).$$
To reach this goal, we first obtain from \eqref{eq-2.5} and \eqref{eq-2.6} that for each $j\in \{1,\cdots,n\}$,
\beq\label{eq-4.2}
&\;\;&\sup_{x\in \mathbb{B}^{n}} \big\{(1-|x|^{2})\big|\nabla (u_{j}\circ \varphi_{y})(x)\big| \big\}=\sup_{x\in \mathbb{B}^{n}} \big\{ \big|\nabla^{h} (u_{j}\circ \varphi_{y})(x)\big| \big\} \\  \nonumber
&=&\sup_{x\in \mathbb{B}^{n}}  \big\{ \left|\nabla^{h} u_{j} \big(\varphi_{y}(x)\big)\right|  \big\}
=\sup_{w\in \mathbb{B}^{n}}  \big\{|\nabla^{h} u_{j}  (w)| \big\} \\  \nonumber
&=&\sup_{w\in \mathbb{B}^{n}} \big\{(1-|w|^{2})\big|\nabla  u_{j} (w)\big| \big\}
<\infty,
\eeq
where in the last inequality, the assumption $u\in \mathcal{B}(\mathbb{B}^{n})$ and \eqref{eq-2.8} are exploited. Then
\begin{eqnarray*}
\big(||u\circ \varphi_{y}||^{\mathcal{B}}\big)^{2}
&=&\sup_{x\in\mathbb{B}^{n} }\left\{(1-|x|^{2})^{2}||D(u\circ \varphi_{y})(x)||^{2}\right\}\\\nonumber
&\leq&\sup_{x\in\mathbb{B}^{n} }\left\{(1-|x|^{2})^{2}\sum_{j=1}^{n}|\nabla (u_{j}\circ \varphi_{y})(x)|^{2}\right\}\;\;\;\;\;\;\;\;\;\;\;\text{(by \eqref{eq-2.8})}\\\nonumber
&\leq&\sum_{j=1}^{n}\left(\sup_{x\in\mathbb{B}^{n} }\left\{(1-|x|^{2})^{2}|\nabla (u_{j}\circ \varphi_{y})(x)|^{2}\right\}\right)\\\nonumber
&=&\sum_{j=1}^{n} \left(\sup_{w\in \mathbb{B}^{n}} \big\{(1-|w|^{2})^{2}\big|\nabla  u_{j} (w) \big|^{2} \big\} \right). \;\;\;\;\;\;\;\;\;\;\;\;\;\;\;\text{(by \eqref{eq-4.2})}
\end{eqnarray*}
Again, by \eqref{eq-2.8} and the assumption $u\in \mathcal{B}(\mathbb{B}^{n})$, we have
\beq\label{eq-4.6}
\;\;\;\big(||u\circ \varphi_{y}||^{\mathcal{B}}\big)^{2}
\leq\sum_{j=1}^{n} \left(\sup_{w\in \mathbb{B}^{n}} \big\{(1-|w|^{2})^{2}|| Du(w)||^{2}\big\} \right)=n\big(||u||^{\mathcal{B}}\big)^{2}<\infty.
\eeq
Hence the claim is proved.\medskip

Now, we have known that $ u\circ \varphi_{y}$ satisfies the conditions in Lemma \ref{lem-4.1}, and so we are ready to finish the proof of the theorem by applying Lemma \ref{lem-4.1}. For $x$, $y\in \mathbb{B}^{n}$, let
$x=\varphi_{y}(z)$. Obviously, $z=\varphi_{y}(x)\in \mathbb{B}^{n}$. Then it follows that
\begin{eqnarray*}
\;\;\;\;\;\;\;\;\;\;\;\;|u(x)-u(y)|
&=& |u\circ \varphi_{y}(z)-u\circ \varphi_{y}(0)|\\\nonumber
&\leq & \frac{||u\circ \varphi_{y}||^{\mathcal{B}} }{2}\log\frac{1+|z|}{1-|z|}\;\;\;\;\;\;\;\;\;\;\;\;\;\;\;\;\;\;\;\;\;(\text{by Lemma 4.1})\\
&\leq &\frac{\sqrt{n}||u||^{\mathcal{B}}}{2} \log\frac{1+|\varphi_{y}(x)|}{1-|\varphi_{y}(x)|} \;\;\;\;\;\;\;\;\;\;\;\;\;\;\;\;\;(\text{by \eqref{eq-4.6}})  \\ \nonumber
&\leq &\frac{\sqrt{n} ||u||_{\mathcal{B}}} {2} \rho(x,y), \;\;\;\;\;\;\;\;\;\;\;\;\;\;\;\;\;\;\;\;\;\;\;\;\;\;\;\;\;(\text{by \eqref{eq-2.1}})
\end{eqnarray*}
as required.\qed

%%%%%%%%%%%%%%%%%%%%%%%%%%%%%%%%%%%%
%%%%%%%%%%%%%%%%%%%%%%%%%%%%%%%%%%%%
\section{Generalized Bloch spaces and integral means}\label{sec-5}
%%%%%%%%%%%%%%%%%%%%%%%%%%%%%%%%%%%%
%%%%%%%%%%%%%%%%%%%%%%%%%%%%%%%%%%%%

The aim of this section is to prove Theorems \ref{thm-1.5} and \ref{thm-1.6}.
Before the proofs, we need some preparation.

\begin{lem}\label{lem-5.1}
Suppose $u\in \mathcal{L}_{p,\omega}\mathcal{B}_{\alpha,\mathrm{a}}^{\beta}(\mathbb{B}^{n})$. Then for $p\in (0,\infty]$ and $r\in [0,1)$,
 $$ M_{p}(r,||Du||)\leq \frac{(\log\mathrm{a})^{\beta} ||u||_{\mathcal{L}_{p,\omega}\mathcal{B}_{\alpha,\mathrm{a}}^{\beta}(\mathbb{B}^{n}) } }{ \omega \big((\log\mathrm{a})^{\beta}\big)\phi_{\alpha, \beta, \mathrm{a}}(r)  } .$$
\end{lem}
\bpf By direct calculations, we see that
$$\phi'_{\alpha, \beta, \mathrm{a}}(r)=(1-r)^{\alpha-1}\left(\log\frac{\mathrm{a}}{1-r}\right)^{\beta-1}\left( \beta-\alpha \log\frac{\mathrm{a}}{1-r}\right).$$
It follows from the assumptions in the lemma that $\phi'_{\alpha, \beta, \mathrm{a}}(r)\leq0$. Hence,
$\phi_{\alpha, \beta, \mathrm{a}}(r)$ is non-increasing in $(0,1)$.
Then the assumption ``$\frac{\omega(t)}{t}$ being non-increasing in $(0,1)$" implies that $ \frac{\phi_{\alpha, \beta, \mathrm{a}}(r)}{\omega(\phi_{\alpha, \beta, \mathrm{a}}(r))}$ is also non-increasing in (0,1).
Therefore,
$$
M_{p}(r,||Du|| )
 \leq   \frac{  ||u||_{\mathcal{L}_{p,\omega}\mathcal{B}_{\alpha,\mathrm{a}}^{\beta}(\mathbb{B}^{n}) }}{ \omega\big(\phi_{\alpha, \beta, \mathrm{a}}(r)\big) }
 \leq   \frac{ ||u||_{\mathcal{L}_{p,\omega}\mathcal{B}_{\alpha,\mathrm{a}}^{\beta}(\mathbb{B}^{n}) }\phi_{\alpha, \beta, \mathrm{a}}(0)}{ \omega \big(\phi_{\alpha, \beta, \mathrm{a}}(0)\big)\phi_{\alpha, \beta, \mathrm{a}}(r) }
 =\frac{(\log\mathrm{a})^{\beta} ||u||_{\mathcal{L}_{p,\omega}\mathcal{B}_{\alpha,\mathrm{a}}^{\beta}(\mathbb{B}^{n}) } }{ \omega \big((\log\mathrm{a})^{\beta}\big)\phi_{\alpha, \beta, \mathrm{a}}(r)  } ,$$
which is what we want.
\epf

%%%%%%%%%%%%%%%%%%%%%%%%%%%%%%%%%%%%
\subsection*{Proof of Theorem \ref{thm-1.5}}
%%%%%%%%%%%%%%%%%%%%%%%%%%%%%%%%%%%%

We divide the proof into two cases according to the values of the parameter $p$.

 \begin{case}\label{case-5.1} $ p\in [1,\infty)$. \end{case}

Let $y=|y|\xi$ and $\gamma(t)=ty$, where $\xi\in \mathbb{S}^{n-1}$ and $t\in[0,1]$. Since
 \be\label{eq-5.1}|u(y)|\leq\ |u(0)|+\int_{\gamma_{[0,y]}}||D u(\gamma)||\cdot |d\gamma| = |u(0)|+\int_{0}^{|y|}||D u(\rho\xi)||  \;d\rho ,\ee
we know from Minkowski's inequality (cf. \cite[Theorem 3.5]{rudin3}) that
$$
M_{p}(|y|,u )=\left( \int_{\mathbb{S}^{n-1}} |u(|y|\xi)|^{p} d\sigma(\xi)\right)^{\frac{1}{p}}\leq |u(0)|+ \left(\int_{\mathbb{S}^{n-1}} \left(\int_{0}^{|y|}||D u(\rho \xi)||  d\rho \right)^{p} d\sigma(\xi)\right)^{\frac{1}{p}}. $$
Further, Minkowski's integral inequality (cf. \cite[A. 1]{stein1970}) gives
$$M_{p}(|y|,u )\leq  |u(0)|+  \int_{0}^{|y|}\left(\int_{\mathbb{S}^{n-1}}||D u(\rho \xi)|| ^{p}d\sigma(\xi)\right)^{\frac{1}{p}}  d\rho.$$

Since it follows from Lemma \ref{lem-5.1} that
\begin{eqnarray*}
 \int_{0}^{|y|}\left(\int_{\mathbb{S}^{n-1}}||D u(\rho \xi)|| ^{p}d\sigma(\xi)\right)^{\frac{1}{p}}  d\rho
&\leq &\frac{(\log\mathrm{a})^{\beta} ||u||_{\mathcal{L}_{p,\omega}\mathcal{B}_{\alpha,\mathrm{a}}^{\beta}(\mathbb{B}^{n}) } }{ \omega \big((\log\mathrm{a})^{\beta}\big) }\int_{0}^{|y|} \frac{1}{ \phi_{\alpha, \beta, \mathrm{a}}(\rho)}   d\rho \\
&=& \frac{(\log\mathrm{a})^{\beta} ||u||_{\mathcal{L}_{p,\omega}\mathcal{B}_{\alpha,\mathrm{a}}^{\beta}(\mathbb{B}^{n}) } }{ \omega \big((\log\mathrm{a})^{\beta}\big) }\int_{0}^{1} \frac{|y|}{ \phi_{\alpha, \beta, \mathrm{a}}(\rho|y|)} d\rho,
 \end{eqnarray*}
we know that for all $y\in \mathbb{B}^{n}$,
\be\label{nd-5}
M_{p}(|y|,u )\leq |u(0)|+\frac{(\log\mathrm{a})^{\beta} ||u||_{\mathcal{L}_{p,\omega}\mathcal{B}_{\alpha,\mathrm{a}}^{\beta}(\mathbb{B}^{n}) } }{ \omega \big((\log\mathrm{a})^{\beta}\big) }\int_{0}^{1} \frac{|y|}{ \phi_{\alpha, \beta, \mathrm{a}}(\rho|y|)} d\rho.
\ee

 \begin{case}\label{case-5.2} $ p=\infty$. \end{case}
For any $y\in \mathbb{B}^{n}$, it follows from the assumption $u\in \mathcal{L}_{p,\omega}\mathcal{B}_{\alpha,\mathrm{a}}^{\beta}(\mathbb{B}^{n})$ and Lemma \ref{lem-5.1} that
 \beq\label{nd-6}
\;\;\;\;\;\;\;\; M_{\infty}(|y|,u) &\leq & |u(0)|+ \int_{0}^{|y|}M_{\infty}(\rho,||D u||) \;d\rho\;\;\;\;\;\;\;\;\;\;\;\;\;\;\;\;\;\;\;\;\text{(by \eqref{eq-5.1})}\\ \nonumber
&\leq & |u(0)|+ \frac{(\log\mathrm{a})^{\beta} ||u||_{\mathcal{L}_{p,\omega}\mathcal{B}_{\alpha,\mathrm{a}}^{\beta}(\mathbb{B}^{n}) } }{ \omega \big((\log\mathrm{a})^{\beta}\big) }\int_{0}^{|y|}\frac{ 1}{\phi_{\alpha, \beta, \mathrm{a}}(\rho)}d\rho\\ \nonumber
&=& |u(0)|+\frac{(\log\mathrm{a})^{\beta} ||u||_{\mathcal{L}_{p,\omega}\mathcal{B}_{\alpha,\mathrm{a}}^{\beta}(\mathbb{B}^{n}) } }{ \omega \big((\log\mathrm{a})^{\beta}\big) }\int_{0}^{1} \frac{ |y|}{ \phi_{\alpha, \beta, \mathrm{a}}(\rho|y|)} d\rho.
 \eeq

We conclude from \eqref{nd-5} and \eqref{nd-6} that the proof of this theorem is complete.
\qed
%
%%%%%%%%%%%%%%%%%%%%%%%%%%%%%%%%%%%%%
%\subsection*{Proof of Corollary \ref{cor-1.1}}
%%%%%%%%%%%%%%%%%%%%%%%%%%%%%%%%%%%%%
%\textcolor[rgb]{1.00,0.00,0.00}{It follows from the} definition that $u\in \mathcal{L}_{p,id}\mathcal{B}_{\alpha,1}^{0}(\mathbb{B}^{n})$ if and only if $$M_{p}(r,||Du||)\leq \frac{ ||u||_{\mathcal{L}_{p,id}\mathcal{B}_{\alpha,1}^{0}(\mathbb{B}^{n}) } }{ (1-r)^{\alpha}} .$$
%Then for $p\in [1,\infty)$ and $|y|\rightarrow 1$, \eqref{eq-5.9} leads to
%\begin{eqnarray*}
%M_{p}(|y|,u )&\leq& |u(0)|+ \int_{0}^{|y|} \frac{ ||u||_{\mathcal{L}_{p,id}\mathcal{B}_{\alpha,1}^{0}(\mathbb{B}^{n})} }{ (1-\rho)^{\alpha}}d\rho
%\leq|u(0)|+ \frac{||u||_{\mathcal{L}_{p,id}\mathcal{B}_{\alpha,1}^{0}(\mathbb{B}^{n})}}{(1-\alpha)(1-|y|)^{\alpha-1}}.
% \end{eqnarray*}
%For $p=\infty$ and $|y|\rightarrow 1$, \eqref{eq-5.1} leads to
%\begin{eqnarray*}
%|u(y)|&\leq& |u(0)|+  \int_{0}^{|y|} \frac{||u||_{\mathcal{L}_{p,id}\mathcal{B}_{\alpha,1}^{0}(\mathbb{B}^{n})} }{ (1-\rho)^{\alpha}}d\rho
%\leq|u(0)|+ \frac{||u||_{\mathcal{L}_{p,id}\mathcal{B}_{\alpha,1}^{0}(\mathbb{B}^{n})}}{(1-\alpha)(1-|y|)^{\alpha-1}}.
% \end{eqnarray*}
%The proof of this theorem is finished.
%\qed

%%%%%%%%%%%%%%%%%%%%%%%%%%%%%%%%%%%%
\subsection*{Proof of Theorem \ref{thm-1.6}}
%%%%%%%%%%%%%%%%%%%%%%%%%%%%%%%%%%%%

We prove this theorem by considering two possibilities according to the values of the parameter $q$.

 \begin{case}\label{case-5.3} $ q=\infty$. \end{case}
To prove the theorem in this case, we need to estimate the operator norm: $||Du(x)||$. We start with some preparation. First, we shall estimate the quantity
$|\nabla u_{j}(x)|^{p}$ in terms of the integral $\int_{\mathbb{B}^{n}\left(0,\frac{1+4|x|}{4+|x|}\right)}| u_{j}(y)|^{p}\;d\nu(y)$. For this, we take $\delta=\frac{1}{4}$ in Lemma \Ref{lem-B}. Then  by \eqref{eq-2.5}, we know that there is a constant $\mu_{09}$ such that
for any $x\in \mathbb{B}^n$ and each $j\in\{1,\cdots,n\}$,
\be\label{eq-5.2}(1-|x|^{2})^{p}|\nabla u_{j}(x)|^{p}\leq \mu_{09}\int_{E(x,\frac{1}{4})}(1-|y|^{2})^{-n} | u_{j}(y)|^{p}\;d\nu(y),\ee
where $\mu_{09}=\mu_{09}\big(\mu_{06}(p,\frac{1}{4})\big)$.
Moreover, by \eqref{eq-2.2}, we obtain that $$ E(x,\frac{1}{4})\subset\mathbb{B}^{n}\left(0,\frac{1+4|x|}{4+|x|}\right).$$
Then we infer from \eqref{eq-2.3} and \eqref{eq-5.2} that
\begin{eqnarray}\label{nd-1}
|\nabla u_{j}(x)|^{p}
&\leq &  \frac{10^{n} \mu_{09} }{3^{n}(1-|x|)^{n+p}}\int_{E(x,\frac{1}{4})} |u_{j}(y)|^{p}\;d\nu(y)\\ \nonumber
&\leq &  \frac{10^{n} \mu_{09} }{3^{n}(1-|x|)^{n+p}}\int_{\mathbb{B}^{n}\left(0,\frac{1+4|x|}{4+|x|}\right)} | u_{j}(y)|^{p} \;d\nu(y) .
\end{eqnarray}

Second, we shall estimate the integral $\int_{\mathbb{B}^{n}\left(0,\frac{1+4|x|}{4+|x|}\right)}| u_{j}(y)|^{p}\;d\nu(y)$.
Since the assumption $$ M_{p}(r,u)=\left( \int_{\mathbb{S}^{n-1}}|u(r\xi)|^{p}d\sigma(\xi) \right)^{\frac{1}{p}} =O\left(\frac{1}{(1-r)^{\alpha}}\right)$$ as $r\rightarrow1^{-}$, together with
the obvious fact $|u(x)|\geq|u_{j}(x)|$ in $\mathbb{B}^{n}$, implies $$ M_{p}(r,u_{j}) \leq M_{p}(r,u),$$
we see that there exists a constant $\mu_{10}$ such that
\be\label{eq-5.5} M_{p}(r,u_{j}) \leq  \frac{ \mu_{10}}{(1-r)^{\alpha}} \ee as $r\rightarrow1^{-}$. It follows from the assumption $\alpha p>1$ that
for each $j\in\{1,\cdots,n\}$,
\begin{eqnarray}\label{nd-2}
&&\int_{\mathbb{B}^{n}\left(0,\frac{1+4|x|}{4+|x|}\right)}| u_{j}(y)|^{p}\;d\nu(y)\\ \nonumber
&= &\int_{0}^{\frac{1+4|x|}{4+|x|} }n\rho^{n-1}\left(\int_{\mathbb{S}^{n-1}}| u_{j}(\rho\xi)|^{p}\;d\sigma(\xi)\right)d\rho \\ \nonumber
&\leq &\mu_{10}^{p}\int_{0}^{\frac{1+4|x|}{4+|x|} } \frac{n\rho^{n-1}}{(1-\rho)^{  \alpha p}}d\rho  \;\;\;\;\;\;\;\;\;\; \;\;\;\;\;\;\;\;\;\;(\text{by \eqref{eq-5.5}})\\ \nonumber
&\leq &\frac{n5^{\alpha p-1}\mu_{10}^{p}}{3^{ \alpha p-1}(\alpha p-1)(1-|x|)^{\alpha p-1}}.\;\;\;\;\;\;\;\;\;\;(\text{since $\alpha p>1$})
\end{eqnarray}

Now, we are ready to get the estimate on $||Du(x)||$. We deduce from \eqref{eq-2.8}, together with the combination of \eqref{nd-1} and \eqref{nd-2}, that
$$||Du(x)||\leq\left(\sum_{j=1}^{n} |\nabla u_{j}(x)|^{2}\right)^{\frac{1}{2}}\leq \sqrt{n} \left(\frac{n  10^{n}5^{\alpha p-1} \mu_{09}\mu_{10}^{p} }{3^{n+\alpha p-1}(\alpha p-1)}\right)^{\frac{1}{p}} \frac{1}{ (1-|x|)^{ \alpha+1+\frac{n-1}{p}}}.$$
Hence
\be\label{mon-1}
M_{\infty}(r,||D u||)=\sup_{\xi\in \mathbb{S}^{n-1}}||Du(r\xi)||=O\left(\frac{ 1}{ (1-r)^{ \alpha+1+\frac{n-1}{p}}}\right)
\ee
as $r\rightarrow 1^{-}$.

 \begin{case}\label{case-5.4} $ q\in(0,\infty)$. \end{case}
Obviously, for all $q\in (0,\infty)$,
 \beq\label{mon-2}
M_{q}(r,||D u||)&= &\left( \int_{\mathbb{S}^{n-1}} || Du(r\xi)||^{q}d\sigma(\xi)\right)^{\frac{1}{q}} \leq M_{\infty}(r,||D u||)\\ \nonumber
&=& O\left(\frac{ 1}{ (1-r)^{ \alpha+1+\frac{n-1}{p}}}\right)
\eeq
as $r\rightarrow 1^{-}$.

Now, we conclude from \eqref{mon-1} and \eqref{mon-2} that for $q\in (0,\infty]$,
$$u\in \mathcal{L}_{q,id} \mathcal{B}_{\alpha+1+\frac{n-1}{p}, \mathrm{a}}^{0}(\mathbb{B}^{n}),$$ and so the proof of this theorem is finished.
  \qed

%%%%%%%%%%%%%%%%%%%%%%%%%%%%%%%%%%%%
\section{Weak uniform boundedness and quasihyperbolic metric}\label{sec-7}

In this section, we shall prove Theorem \ref{thm-1.8}.

%%%%%%%%%%%%%%%%%%%%%%%%%%%%%%%%%%%%
\subsection*{Proof of Theorem \ref{thm-1.8}}
%%%%%%%%%%%%%%%%%%%%%%%%%%%%%%%%%%%%

The sufficiency in the theorem easily follows from the proof of \cite[Theorem 2.8]{matel} or the proof of \cite[Theorem 2]{chenantti2016}.
So we only need to prove the necessity.
Obviously, for every $x\in \mathbb{B}^{n}$ and $y\in \overline{\mathbb{B}}^{n}\left(x,\frac{d_{\mathbb{B}^{n}}(x)}{4}\right)$,
$$d_{\mathbb{B}^{n}}(y)\geq d_{\mathbb{B}^{n}}(x)-\frac{1}{4}d_{\mathbb{B}^{n}}(x)= \frac{3}{4}d_{\mathbb{B}^{n}}(x).$$ Hence
$$|x-y|\leq \frac{d_{\mathbb{B}^{n}}(x)  }{4}\leq \min\left\{ \frac{d_{\mathbb{B}^{n}}(x)}{2},  \frac{ d_{\mathbb{B}^{n}}(y)}{2}\right\},$$ and
so $$r_{ \mathbb{B}^{n}}(x,y)\leq \frac{1}{2}.$$
Then the assumption ``$u$ satisfying the weak uniform boundedness property" implies that
$$|u(y)-u(x)|\leq \mu_{02}d_{u(\mathbb{B}^{n})}\big(u(x)\big),$$
where $\mu_{02}$ is the constant from \eqref{eq-1.4}.
From Lemma \ref{lem-3.1}, we deduce that
\be\label{eq-7.1}
||D u(x)||\leq
\sqrt{n}\mu_{02}\mu_{07} |\mathbb{B}^n|\frac{5^{n}\cdot d_{u(\mathbb{B}^{n})}\big(u(x)\big)}{8^{n}\cdot d_{\mathbb{B}^{n}}(x) }.
\ee

 Since $|du(z)|\leq ||D u(z)|| \cdot|dz|,$ by taking
 $$\mu_3=\frac{5^{n} }{8^{n} }\sqrt{n}\mu_{02}\mu_{07}|\mathbb{B}^n|,$$
 it follows from \eqref{eq-7.1} that for all $x,$ $y\in \mathbb{B}^{n}$,
\begin{eqnarray*}
k_{u\left(\mathbb{B}^{n}\right)}(u(x),u(y))
&\leq&\inf_{\gamma\in\Gamma_{xy}(\mathbb{B}^{n})}\int_{\gamma}\frac{1}{d_{u(\mathbb{B}^{n})}\big(u(z)\big)}|du(z)|\\
&\leq& \inf_{\gamma\in\Gamma_{xy}(\mathbb{B}^{n})}\int_{ \gamma}\frac{||D u(z)||}{d_{u(\mathbb{B}^{n})}\big(u(z)\big)}|dz| \\
&\leq&\inf_{\gamma\in\Gamma_{xy}(\mathbb{B}^{n})}\int_{ \gamma}\frac{\mu_{3}}{d_{\mathbb{B}^{n}}(z) }|dz|\\
&= &\mu_{3} k_{\mathbb{B}^{n}}(x,y),
\end{eqnarray*}
and so Theorem \ref{thm-1.8} is proved. \qed

%\bigskip
%\noindent {\bf Acknowledgements:}
%The research was partly supported by NSFs of China (No. 11571216 and No. 11671127).

\end{document}